\documentclass[12pt]{amsart}
\usepackage[margin=3cm, marginpar=3cm]{geometry}

\usepackage{amssymb}
\usepackage{braket}
\usepackage{mathrsfs}
\usepackage{ifthen}
\usepackage{graphicx}

\usepackage{comment} 
\usepackage{mleftright}
\usepackage[pagebackref,hypertexnames=false]{hyperref} 
\usepackage{cleveref}
 \usepackage[all]{xy}
 \usepackage{amscd}
 \usepackage{color}
 \usepackage{enumitem}
\newlist{steps}{enumerate}{1}
\setlist[steps, 1]{label = Step \arabic*:}
\usepackage[fontsize=10]{scrextend}
\usepackage[pagebackref,hypertexnames=false]{hyperref} 
\usepackage[alphabetic,backrefs,msc-links]{amsrefs}

\makeatletter
\DeclareRobustCommand\widecheck[1]{{\mathpalette\@widecheck{#1}}}
\def\@widecheck#1#2{%
   \setbox\z@\hbox{\m@th$#1#2$}%
   \setbox\tw@\hbox{\m@th$#1%
      \widehat{%
         \vrule\@width\z@\@height\ht\z@
         \vrule\@height\z@\@width\wd\z@}$}%
   \dp\tw@-\ht\z@
   \@tempdima\ht\z@ \advance\@tempdima2\ht\tw@ \divide\@tempdima\thr@@
   \setbox\tw@\hbox{%
      \raise\@tempdima\hbox{\scalebox{1}[-1]{\lower\@tempdima\box\tw@}}}%
   {\ooalign{\box\tw@ \cr \box\z@}}}
\makeatother

\theoremstyle{plain}
\newtheorem{thm}{Theorem}[section]
\crefname{thm}{Theorem}{Theorems}
\Crefname{thm}{Theorem}{Theorems}

\crefname{prop}{Proposition}{Propositions}
\Crefname{prop}{Proposition}{Propositions}

\crefname{lem}{Lemma}{Lemmas}
\Crefname{lem}{Lemma}{Lemmas}

\crefname{cor}{Corollary}{Corollaries}
\Crefname{cor}{Corollary}{Corollaries}

\crefname{claim}{Claim}{Claims}
\Crefname{claim}{Claim}{Claims}

\crefname{property}{Property}{Properties}
\Crefname{property}{Property}{Properties}

\crefname{problem}{Problem}{Problems}
\Crefname{problem}{Problem}{Problems}

\crefname{ques}{Question}{Questions}
\Crefname{ques}{Question}{Questions}

\theoremstyle{definition}

\crefname{defn}{Definition}{Definitions}
\Crefname{defn}{Definition}{Definitions}

\crefname{notation}{Notation}{Notations}
\Crefname{notation}{Notation}{Notations}

\crefname{convention}{Convention}{Conventions}
\Crefname{convention}{Convention}{Conventions}

\crefname{cond}{Condition}{Conditions}
\Crefname{cond}{Condition}{Conditions}

\crefname{assum}{Assumption}{Assumptions}
\Crefname{assum}{Assumption}{Assumptions}

\crefname{conj}{Conjecture}{Conjectures}
\Crefname{conj}{Conjecture}{Conjectures}

\theoremstyle{remark}
\newtheorem{rem}[thm]{Remark}
\crefname{rem}{Remark}{Remarks}
\Crefname{rem}{Remark}{Remarks}

\crefname{ex}{Example}{Examples}
\Crefname{ex}{Example}{Examples}

\crefname{section}{Section}{Sections}
\Crefname{section}{Section}{Sections}
\crefname{subsection}{Subsection}{Subsections}
\Crefname{subsection}{Subsection}{Subsections}
\crefname{figure}{Figure}{Figures}
\Crefname{figure}{Figure}{Figures}

\newcommand{\Z}{\mathbb{Z}}

\newcommand{\C}{\mathbb{C}}

\newcommand{\ctext}[1]{\raise0.2ex\hbox{\textcircled{\scriptsize{#1}}}}

\newcommand{\mbar}[1]{{\ooalign{\hfil#1\hfil\crcr\raise.167ex\hbox{--}}}}

\title{On the connected sums of the $(2,1)$-cable of the figure eight knot}
\author{Yoshihiro Fukumoto and Masaki Taniguchi }

\begin{document}

\maketitle

\begin{abstract}
   We show that the 3-fold (resp. 6-fold) connected sum of the \((2,1)\)-cable of the figure-eight knot cannot bound a smooth null-homologous disk in a punctured \(S^2 \times S^2\) (resp. in a punctured \(\#_2 S^2 \times S^2\)). This result is obtained using a real version of the \(10/8\)-inequality established by Konno, Miyazawa, and Taniguchi.
\end{abstract}

\section{Introduction}
The \((2,1)\)-cable of the figure-eight knot, denoted by \((4_1)_{(2,1)}\), has long been known to be algebraically slice and strongly rationally slice~\cites{kawauchi1, Cha:2007-1, Kim-Wu:2018-1}. It has also been considered as a potential counterexample to the slice-ribbon conjecture~\cite[Problem 25]{Fox:1962-1}. However, most known concordance invariants failed to detect the non-sliceness of this knot, leaving its smooth sliceness an open question.

Recently, Dai, Kang, Mallick, Park, and Stoffregen showed in~\cite[Theorem 1.1]{DKMPS:2022-1} that \((4_1)_{(2,1)}\) is not smoothly slice by employing involutive Heegaard Floer theory. 
(See ~\cite[Theorem 2.1]{ACMPS:2023-1}, ~\cite[Corollary 1.20]{KMT23f} and ~\cite{KPT24} for alternative proofs.) Furthermore, they demonstrated that the \(n\)-fold connected sum $\#_n(4_1)_{(2,1)}$ is not smoothly slice for any \(n \neq 0\), thereby strengthening the result.

 In this article, we shall prove: 
\begin{thm}\label{main thm}
 The $3$-fold (resp. the $6$-fold) connected sum of $(2,1)$-cable of the figure eight knot does not bound smooth null-homologous disk in a punctured $ S^2\times S^2$ (resp. in a punctured $\#_2S^2\times S^2$). 
 \end{thm}

It is shown in \cite[page 84, Corollary~5.11]{Ka87} and \cite{Sch10} that for any knot \(K\) in \(S^3\) with vanishing Arf invariant \(\operatorname{Arf}(K) = 0\), the knot \(K\) bounds a smooth null-homologous disk in a punctured \(\#_l S^2 \times S^2\) for some \(l \geq 0\). The minimum of such \(l\) is denoted by \(\operatorname{sn}(K)\) and is called the \emph{stabilizing number} of the knot \(K\). \Cref{main thm} can be restated as follows: 
\[
\operatorname{sn}(\#_3(4_1)_{(2,1)}) \geq 2 \quad \text{and} \quad \operatorname{sn}(\#_6(4_1)_{(2,1)}) \geq 3.
\]

\begin{rem}
    Our techniuqe relies on the existence of a smooth concordance from the figure-eight knot to the unknot in a  twice-punctured $2\mathbb{CP}^2$, denoted by $X$, that represents $(1,3)$ in $H_2(X, \partial X; \mathbb{Z}) \cong \mathbb{Z} \oplus \mathbb{Z}$, as proved by Aceto, Castro, Miller, Park, and Stipsicz in~\cite{ACMPS:2023-1}. So, our method applies to all knots that permit such a concordance to a smoothly slice knot, which is the case for \cite[Theorem 2.3]{ACMPS:2023-1} as well.
    Therefore, one can see: 
    Let $K$ be a knot, and let $K_{(2,1)}$ denote the $(2,1)$-cable of $K$. Suppose that $K$ can be transformed into a slice knot by applying full negative twists along two disjoint disks, where one intersects $K$ algebraically once and the other intersects it algebraically three times. Then we have 
    \[
\operatorname{sn}(\#_3K_{(2,1)}) \geq 2 \quad \text{and} \quad \operatorname{sn}(\#_6K_{(2,1)}) \geq 3.
\]
There are infinitely many knots that satisfy the assumptions, see \cite[Remark 2]{ACMPS:2023-1}.  
\end{rem}

The proof utilizes a real version of the Manolescu's relative \(10/8\)-inequality \cite{Ma14}  generalizing Furuta's 10/8 inequality \cite{Fu01}, which gives constraints on surfaces smoothly embedded in 4-manifolds, as established in \cite[Theorem 1.3(v)]{KMT21}. 
Alternatively, the real 10/8-inequality in \cite{KMT21} is regarded as a relative version of a variant 10/8 inequality proven by Kato \cite{Kat22}. For the background of real Seiberg--Witten theory, see \cite{TW09, Na13, Nak15, Ka22, KMT21, Ji22, KMT:2023, Mi23, Li23}.

This real 10/8 inequality \cite[Theorem 1.3(v)]{KMT21} for surfaces is derived by considering double-branched covers and applying the real \(10/8\)-inequality \cite[Theorem 1.1(iv)]{KMT21} for smooth involutions. Furthermore, our approach yields the following result for involutions:

\begin{thm}\label{involution}
    Let $Y$ be the double-branched cover along $(4_1)_{(2,1)}$.  Denote by $\tau$ the covering involution on $ Y$.  Then, for any smooth compact $\Z_2$-homology $4$-ball $X$ bounded by $\#_3Y$ (resp. $\#_6Y$), the action $\#_3\tau$ (resp. $\#_6\tau$) does not extend to $X\#S^2\times S^2$ (resp. $X\# \#_2 S^2\times S^2$) as a smooth $\Z_2$-action.  
\end{thm}
The double branched cover along $(4_1)_{(2,1)}$ is known to be diffeomorphic to $S^3_1(4_1 \# 4_1^r)$, where $K^r$ denotes $K$ with the reversed orientation. Since $4_1 \# 4_1^r$ is smoothly slice, $Y$ bounds a smooth compact contractible $4$-manifold. 
Note that \(\#_nY\) with \(\#_n\tau\) is shown to be a {\it strong cork} in~\cite[Theorem 1.1]{DKMPS:2022-1} for $n \neq 0$. From this perspective, \cref{involution} can be compared with the existence of strong corks that survive after one stabilization, as established by Kang in \cite{Kang22}, again based on involutive Heegaard Floer theory. \cref{involution} suggests that \(\#_6 Y\) with \(\#_6 \tau\) provides a candidate for a strong cork that survives after two stabilizations. Once the existence of strong corks surviving after two stabilizations is established, from Akbulut--Ruberman's technique \cite{AR16}, one can conclude the existence of compact contractible exotic 4-manifold pair surviving after two stabilizations.

\subsection*{Acknowledgment} 
The first author was partially supported by JSPS KAKENHI Grant Number 22K03322. 
The second author was partially supported
by JSPS KAKENHI Grant Number 22K13921.

\section{Proof of the main theorems}

As we mentioned in the introduction, the main tool to prove is a version of real 10/8 inequality \cite[Theorem 3.41]{KMT21}, which is, slightly stronger than the usual real 10/8 inequality \cite[Theorem 1.1]{KMT21}. Also, we shall combine it with the surface cobordism constructed in \cite{ACMPS:2023-1}. 
\begin{proof}[Proof of \cref{main thm}]
We start with a smooth surface cobordism $S$ of genus $0$ in $W:= [0,1]\times S^3 \# (2\C P^2)$ from $(4_1)_{2, 1}$ to $T(2, -19)$ whose homology class is 
\[
(2 , 6 ) \in H_2(W; \Z) = H_2( 2\C P^2; \Z)= \Z \oplus \Z, 
\]
which is constructed in \cite{ACMPS:2023-1}. Here $T(2, -19)$ denotes the torus knot of type $(2,-19)$.

Let me be a positive integer. 
We consider the $m$-fold boundary connected sums 
\[
S_m := \natural_m S \subset [0,1]\times S^3 \# 2m \C P^2 = W_m
\]
of the surface $S$ along certain paths so that $\natural_m S$ gives a surface cobordism from $\#_m (4_1)_{2, 1}$ to $\#_m T(2, -19)$. We take the double-branched cover $\Sigma_2(S_m)$ along $S_m$ of $[0,1]\times S^3 \# (2m \C P^2)$ which is regarded as a $\Z_2$-equivariant cobordism 
\[
\partial \Sigma_2(S_m) = - \Sigma_2(\#_m (4_1)_{2, 1})  \cup - \#_m\Sigma(2, 2, 19) ,
\]
where $\Sigma(2, 2, 19)$ is the Brieskorn $3$-sphere of type $(2, 2, 19)$.  Since $PD(w_2(W)) \equiv [S_m] \operatorname{mod} 2$ holds, $\Sigma_2(S_m)$ is spin and one can check it has unique spin structure. 
One can see that 
\[
b^+(\Sigma_2(S_m)) - b^+(W_m)  = m (b^+(\Sigma_2(S)) - b^+(W))=m  \text{ and }\sigma( \Sigma_2(S_m)) = m \sigma(\Sigma_2(S))=  2m.
\]
For the computations of these quantities, see \cite[Lemma 4.2]{KMT21} for example. See \cite[Section 3.4]{KPT24} for related computations. 

Suppose $\#_m (4_1)_{(2,1)}$ bounds a smoothly embedded null-homologous disk in punctured $\#_lS^2\times S^2$ for some positive integer $l$. We take a double-branched cover along the disk and get a $\Z_2$-equivariant $4$-manifold $X$ bounded by $\Sigma_2(\#_m (4_1)_{(2,1)})$. Note that $X$ is $\Z_2$-homology $D^4$, hence it is spin.

Consider the spin $\Z_2$-equivariant 4-manifold: 
\[
Z_m := X \cup_{\Sigma_2(\#_m (4_1)_{(2,1)})} \Sigma_2(S_m). 
\]
Note that $Z_m$ has unique spin structure. 
Again from \cite[Lemma 4.2]{KMT21}, one can check 
\[
b^+ (Z_m) - b^+(Z_m/\Z_2) = m + l \text{ and }\sigma(\Sigma_2(Z_m)) = 2 m .
\]

We apply the real 10/8 inequality \cite[Theorem 3.41]{KMT21} to $Z_m$ with the involution and the odd spin structure and, under suitable assumption, obtain 
\begin{align}\label{strong108}
  - \frac{\sigma (\Sigma_2(Z_m))}{16} + A (N) \leq b^+(Z_m) - b^+(Z_m/\Z_2) + \kappa_R ( -\#_m \Sigma_2(2,2,19)), 
\end{align}
where $\kappa_R$ denotes the real kappa invariant of rational spin homology $3$-spheres with odd involutions. 
For $N \in \Z$,  
$A(N) \in \{1,2,3\}$ is defined by
\begin{align}\label{A-definition}
A(N) = 
\begin{cases}
&1 , \quad N=0,2 \operatorname{mod} 8\\
&2, \quad N=1,3,4,5,7 \operatorname{mod} 8\\
&3, \quad N=6 \operatorname{mod} 8.
\end{cases}    
\end{align}
Here 
\[
N:=
-\frac{\sigma(\Sigma_2(Z_m))}{16} -\kappa_R ( -\#_m \Sigma_2(2,2,19)).
\]
Also, in order to get \eqref{strong108}, we need to check
\[
N\geq 2 \quad \text{and} \quad 
b^+(Z_m) - b^+(Z_m/\Z_2) \geq 1,
\] 
which are the assumptions of \cite[Theorem 3.41]{KMT21}.

Since $\Sigma_2(2,2,19)$ is lens space, from \cite[Example 3.56]{KMT21} and the additivity \cite[Corollary 3.51]{KMT21} of the real kappa invariant, we can see 
\[
\kappa_R ( -\#_m \Sigma_2(2,2,19)) =  \frac{1}{16} \sigma(\#_m T(2,19))  = -\frac{9}{8}m ,  
\]
where $\sigma (K)$ denotes the knot signature with the convention $\sigma(T(2,3)) =-2$.
So, in our situation, we have $N = m$ and $b^+(Z_m) - b^+(Z_m/\Z_2) = m+ l$.   Assuming that $N=m \geq 2$, we have 
\[
\begin{cases}
&1 \leq l , \quad m=0,2 \operatorname{mod} 8\\
&2 \leq l , \quad m=1,3,4,5,7 \operatorname{mod} 8\\
&3 \leq l , \quad m=6 \operatorname{mod} 8.
\end{cases}      
\]

Putting $m = 3$ and $6$, we get the desired results. 
\end{proof}

\begin{rem}
Note that this arguemnt also ensures $ \#_m (4_1)_{(2,1)}$ is not slice for $m \geq 2$. Let us point out that \cite[Theorem 3.41]{KMT21} is not enough to directly prove $ (4_1)_{(2,1)}$ is not slice although non-sliceness of $ \#_2 (4_1)_{(2,1)}$ implies that of $ (4_1)_{(2,1)}$. 
Using a technique to see the inequality proven in \cite[Theorem 1.3]{Kon24}, one can still check $ (4_1)_{(2,1)}$ is not slice directly from the real 10/8 inequality. 
\end{rem}

\begin{proof}[Proof of \cref{involution}]
The proof is similar to that of \cref{main thm}. Let $X'$ be a smooth $\Z_2$-homology $D^4$ bounded by $\#_m Y$ with an extension of $\#_m\tau$. Then we consider a spin $4$-manifold
\[
Z'_m := X' \cup_{\Sigma_2\left(\#_m (4_1)_{(2,1)}\right)} \Sigma_2(S_m)
\]
with an odd involution instead of $Z_m$. Then the latter proof is completely the same.  
\end{proof}

\bibliographystyle{alpha}
\bibliography{tex}

\end{document}